arXiv:math.CT/0101154 18 Jan 2001

# On the monad of proper factorisation systems in categories (*)

Marco Grandis

*Dipartimento di Matematica, Università di Genova, Via Dodecaneso 35, 16146-Genova, Italy*

*E-mail address:* grandis@dima.unige.it (Marco Grandis).

**Abstract.** It is known that factorisation systems in categories can be viewed as unitary pseudo algebras for the monad $(-)^{\mathbf{2}}$, in **Cat**. We show in this note that an analogous fact holds for *proper* (i.e., *epi-mono*) factorisation systems and a suitable quotient of the former monad, deriving from a construct introduced by P. Freyd for stable homotopy. Structural similarities of the previous monad with the path endofunctor of topological spaces are considered.

*MSC*: 18A32; 18C15.

*Key Words*: Factorisation systems, 2-monads, Eilenberg-Moore algebras, pseudo algebras.

## Introduction

For a category $\mathbf{X}$, the category of morphisms $\mathrm{P}\mathbf{X} = \mathbf{X}^{\mathbf{2}}$ has a natural factorisation system. So equipped, it is the free category with factorisation system, on $\mathbf{X}$.

This system induces a *proper*, or *epi-mono*, factorisation system on a quotient $\mathrm{Fr}\mathbf{X} = \mathbf{X}^{\mathbf{2}}/\mathrm{R}$ [G3], the free category with epi-mono factorisation system on $\mathbf{X}$ (the *epi-mono completion*), that generalises the Freyd embedding of the stable homotopy category of spaces in an abelian category [Fr]. "Weak subobjects" in $\mathbf{X}$, of interest for homotopy categories, correspond to ordinary subobjects in $\mathrm{Fr}\mathbf{X}$; other results in [G3] concern various properties of $\mathrm{Fr}\mathbf{X}$ that derive from weak (co)limits of $\mathbf{X}$.

Now, the "path" endofunctor $\mathrm{P} = (-)^{\mathbf{2}}$ of **Cat** has an obvious 2-monad structure (with diagonal multiplication), linked to the universal property recalled above (a pseudo adjunction); it is known, since some hints in Coppey [Co] and a full proof in Korostenski - Tholen [KT], that its pseudo algebras correspond to the factorisation systems of $\mathbf{X}$. Similarly, as stated without proof in [G3], the pseudo algebras for the induced 2-monad on $\mathrm{Fr}\mathbf{X}$ correspond to *proper* factorisation systems of $\mathbf{X}$; more precisely, we prove here, in Theorem 4 (ii), that there is a canonical bijection between *proper* factorisation systems in $\mathbf{X}$ and pseudo isomorphism classes of unitary pseudo $\mathrm{Fr}$-algebras on $\mathbf{X}$. Similar, simpler relations hold in the *strict* case: *strict factorisation systems are monadic on categories, as well as the proper such*. Structural similarities of $\mathrm{P}$ with the topological path functor $\mathrm{PX} = \mathrm{X}^{[0,1]}$ are discussed at the end (Section 5).

We shall use the same notation of [G3]. For factorisation systems, one can see Freyd - Kelly [FK], Carboni - Janelidze - Kelly - Paré [CJKP], and their references; the strict version is much less used: see [G3] and Rosebrugh-Wood [RW]. Lax $\mathrm{P}$-algebras are studied in [RT]. General lax and pseudo algebras can be found in Street [St].

---

(*) Work partially supported by MURST research projects.

The author acknowledges with pleasure a suggestion of F.W. Lawvere, at the origin of this note.

**1. The factorisation monad.** Let **X** be any category and $\mathbf{X}^2$ its category of morphisms. An object of the latter is an **X**-map x: X' → X", which we *may* write as $\hat{x}$ when it is viewed as an object of $\mathbf{X}^2$; a morphism f = (f', f"): $\hat{x} \to \hat{y}$ is a commutative square of **X**, as in the left diagram

$$(1)\qquad \begin{array}{ccc} X' & \xrightarrow{f'} & Y' \\ {\scriptstyle x}\downarrow & & \downarrow{\scriptstyle y} \\ X'' & \xrightarrow{f''} & Y'' \end{array} \qquad\qquad \begin{array}{ccccc} X' & = & X' & \xrightarrow{f'} & Y' \\ {\scriptstyle x}\downarrow & & \downarrow{\scriptstyle \bar{f}} & & \downarrow{\scriptstyle y} \\ X'' & \xrightarrow{f''} & Y'' & = & Y'' \end{array}$$

and the composition is obvious. The *strict factorisation* of f, shown in the right diagram, is f = (f', 1).(1, f"); note that its middle object is the *diagonal* $\bar{f}$ = f"x = yf' of the square f.

Thus, $\mathbf{X}^2$ has a canonical factorisation system (*fs* for short), where the map f = (f', f") is in E (resp. in M) iff f' (resp. f") is an isomorphism. This system contains a canonical *strict factorisation system,* where (f', f") is in $E_0$ (resp. in $M_0$) iff f' (resp. f") is an identity. (As in [G3, 2.1], this means that: (i) $E_0$, $M_0$ are subcategories containing all the identities; (ii) every map u has a *strictly unique* factorisation u = me with $e \in E_0$, $m \in M_0$. A strict fs $(E_0, M_0)$ is not a fs, of course; but, there is a unique fs (E, M) containing the former, where u = me is in E iff m is iso, and dually. Two strict systems are said to be *equivalent* if they span the same fs.)

The full embedding that identifies the object X of **X**, with $\hat{1}_X$

$$(2)\qquad \eta\mathbf{X}: \mathbf{X} \to \mathbf{X}^2, \qquad\qquad (f: X \to Y) \mapsto (f, f): \hat{1}_X \to \hat{1}_Y,$$

makes $\mathbf{X}^2$ the *free category with factorisation system* on **X**, in the "ordinary" sense (as well as in a *strict* sense): for every functor F: **X** → **A** with values in a category with fs (resp. *strict* fs), there is an extension G: $\mathbf{X}^2 \to \mathbf{A}$ that preserves factorisations (resp. *strict* factorisations), determined up to a unique functorial isomorphism (resp. *uniquely* determined): $G(\hat{x}) = \mathrm{Im}_{\mathbf{A}}(Fx)$. The (obvious) proof is based on the canonical factorisation of $\eta\mathbf{X}(x) = (x, 1).(1, x)$: $\hat{1}_{X'} \to \hat{1}_{X''}$ in $\mathbf{X}^2$

$$(3)\qquad \begin{array}{ccccc} X' & = & X' & \xrightarrow{x} & X'' \\ {\scriptstyle 1}\downarrow & & \downarrow{\scriptstyle x} & & \downarrow{\scriptstyle 1} \\ X' & \xrightarrow{x} & X'' & = & X'' \end{array} \qquad\qquad X' \xrightarrow{(1,\, x)} \hat{x} \xrightarrow{(x,\, 1)} X''.$$

One might now expect that "factorisation systems be monadic on categories", but this is only true in a relaxed 2-dimensional sense.

First, by the strict universal property, the forgetful 2-functor $U_0$: $\mathrm{Fs}_0$**Cat** → **Cat** (of categories with *strict* fs) has a left 2-adjoint $F_0(\mathbf{X}) = (\mathbf{X}^2; E_0, M_0)$, and we shall see that $U_0$ is indeed 2-monadic: the comparison 2-functor $K_0$: $\mathrm{Fs}_0$**Cat** → P-**Alg** establishes an isomorphism of $\mathrm{Fs}_0$**Cat** with the 2-category of algebras of the associated 2-monad, $P = U_0F_0$: **Cat** → **Cat**, $P(\mathbf{X}) = \mathbf{X}^2$.

Secondly, by the "relaxed" universal property, the forgetful 2-functor U: Fs**Cat** → **Cat** (of categories with fs) acquires a left *pseudo* adjoint 2-functor $F(\mathbf{X}) = (\mathbf{X}^2; E, M)$: the unit η: 1 → UF is 2-natural, but the counit is *pseudo* natural and "ill-controlled", each component $\varepsilon\mathbf{A}$: $(|\mathbf{A}|^2; E, M) \to \mathbf{A}$ depending on a choice of images in **A**; the triangle conditions are – rather – invertible 2-cells. This would give an ill-determined *pseudo* monad structure on $P = UF = U_0F_0$, *isomorphic* to the previous 2-monad; we will therefore settle on the latter and "by-pass" the pseudo adjunction.

In fact, the structure of the category $\mathbf{2} = \{0 \to 1\}$ as a *diagonal* comonoid (with $e: \mathbf{2} \to \mathbf{1}$, $d: \mathbf{2} \to \mathbf{2}{\times}\mathbf{2}$) produces a *diagonal* monad on the endofunctor $P = (-)^\mathbf{2}$ of **Cat**, precisely the one we are interested in. The unit $\eta\mathbf{X} = \mathbf{X}^e: \mathbf{X} \to \mathbf{X}^\mathbf{2}$ is the canonical embedding considered above, $\eta\mathbf{X}(X) = \hat{1}_X$. The multiplication $\mu\mathbf{X} = \mathbf{X}^d: P^2\mathbf{X} \to P\mathbf{X}$ is a "diagonal functor" defined on $P^2\mathbf{X} = \mathbf{X}^{\mathbf{2}\times\mathbf{2}}$:

- an object of $P^2\mathbf{X}$ is a morphism $\xi_0 = (a_0, b_0): x_0 \to y_0$ of $P\mathbf{X}$, and a commutative square in $\mathbf{X}$ (the front square of the diagram below); $\mu\mathbf{X}(\xi_0) = d_0 = b_0x_0 = y_0a_0$ is the diagonal of this square;

- a morphism of $P^2\mathbf{X}$ is a commutative square $\Xi$ of $P\mathbf{X}$, and a commutative cube in $\mathbf{X}$; $\mu\mathbf{X}(\Xi)$ is a diagonal square of the cube

(4)

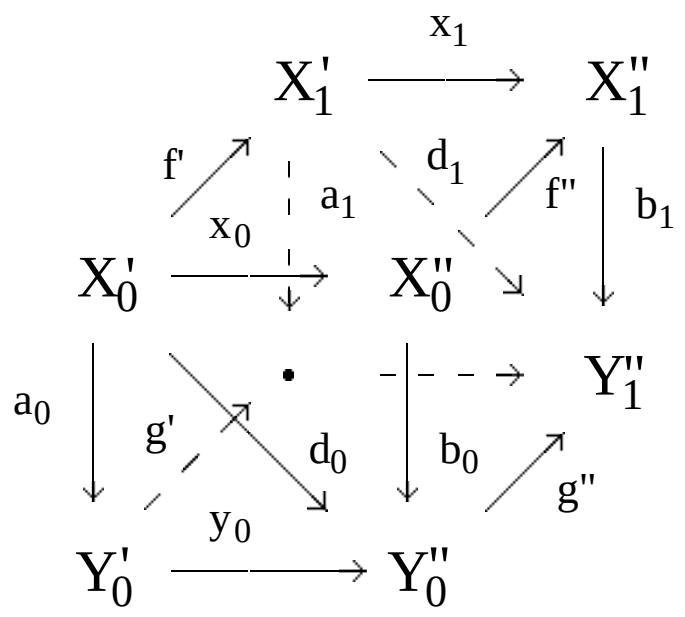

$\Xi = ((f', g'), (f", g")): \xi_0 \to \xi_1,$

$\xi_i = (a_i, b_i): x_i \to y_i,$

$\mu\mathbf{X}(\Xi) = (f', g"): d_0 \to d_1;$

$\mu$ coincides with the multiplication coming from the strict adjunction, $U_0\varepsilon_0F_0 : P^2 \to P$ (and *would* also coincide with the pseudo multiplication $U\varepsilon F$, *if* one might control the choice of images in $F\mathbf{X}$ by its strict fs).

$P$ will also be called the *factorisation monad* on **Cat**, while a $P$-algebra $(\mathbf{X}, t)$ will also be called a *factorisation algebra*; it consists of a functor $t: \mathbf{X}^\mathbf{2} \to \mathbf{X}$ such that $t.\eta\mathbf{X} = 1_\mathbf{X}$, $t.Pt = t.\mu\mathbf{X}$.

**2. The proper factorisation monad.** Consider now the quotient $Fr\mathbf{X} = \mathbf{X}^\mathbf{2}/R$, modulo the "Freyd congruence" [Fr]: two parallel $\mathbf{X}^\mathbf{2}$-morphisms $f = (f', f"): x \to y$ and $g = (g', g"): x \to y$ are R-equivalent whenever their diagonals $\overline{f}$, $\overline{g}$ coincide (cf. 1.1); the morphism of $Fr\mathbf{X}$ represented by $f$ will be written as $[f]$ or $[f', f"]$. As a crucial effect of this congruence, if $f'$ is epi (resp. $f"$ is mono) in $\mathbf{X}$, so is $[f]$ in $Fr\mathbf{X}$.

As in [G3], a *canonical epi* (resp. *mono*) of $Fr\mathbf{X}$ will be a morphism which can be represented as $[1, f"]$ (resp. $[f', 1]$). Every map $[f]$ has a precise *canonical factorisation* $[f] = [f', 1].[1, f"]$, formed of a canonical epi and a canonical mono (both their diagonals being $\overline{f}$). $Fr\mathbf{X}$ has thus a *proper* strict fs $(E_0, M_0)$, which spans a (proper) fs $(E, M)$: the map $[f]: x \to y$ belongs to $E$ iff there is some $u: Y' \to X'$ such that $yf'u = y$ ($y$ *sees* $f'$ *as a split epi*).

The full embedding $\eta'\mathbf{X} = p.\eta\mathbf{X}: \mathbf{X} \to Fr\mathbf{X}$ takes $f: X \to Y$ to $[f, f]: \hat{1}_X \to \hat{1}_Y$; $Fr\mathbf{X}$ is thus *the free category with proper factorisation system* on $\mathbf{X}$ [G3, 2.3], called the *Freyd completion*, or *epi-mono completion* of $\mathbf{X}$. The 2-monad structure of $Fr$, induced by the one of $P$ (by-passing again a pseudo adjunction $F' \dashv U'$), will be called the *proper factorisation monad* on **Cat**. The unit is $\eta'$. For the multiplication $\mu'\mathbf{X}: Fr^2\mathbf{X} \to Fr\mathbf{X}$, note that now

- an object of $Fr^2\mathbf{X}$ is a morphism of $Fr\mathbf{X}$, $\xi_0 = [a_0, b_0]: x_0 \to y_0$,

- a morphism of $Fr^2\mathbf{X}$ is an equivalence class $\Xi$ of commutative squares of $Fr\mathbf{X}$

(1) $\Xi = [[f', g'], [f'', g'']]\colon (\xi_0\colon x_0 \to y_0) \to (\xi_1\colon x_1 \to y_1)$,

and we have

(2) $\mu'\mathbf{X}(\xi_0) = d_0, \qquad\qquad \mu'\mathbf{X}(\Xi) = [f', g'']\colon d_0 \to d_1$;

in fact, the class $[f', g'']\colon d_0 \to d_1$ is well defined, since its diagonal $g''d_0 = g''b_0x_0$ only depends on the class $[f'', g'']$ and the object $x_0$. The projection $p$ is thus a strict morphism of monads $(\mathrm{P}, \eta, \mu) \to (\mathrm{F}r, \eta', \mu')$, as shown in the left diagram below (with $p_2 = \mathrm{F}r(p).p\mathrm{P} = p\mathrm{F}r.\mathrm{P}(p)$)

$$(3) \qquad \begin{array}{ccccc} \mathbf{X} & \xrightarrow{\eta} & \mathrm{P}\mathbf{X} & \xleftarrow{\mu} & \mathrm{P}^2\mathbf{X} \\ \| & & \downarrow p & & \downarrow p_2 \\ \mathbf{X} & \xrightarrow{\eta'} & \mathrm{Fr}\mathbf{X} & \xleftarrow{\mu'} & \mathrm{Fr}^2\mathbf{X} \end{array} \qquad\qquad \begin{array}{ccccc} \mathbf{X} & \xrightarrow{\eta'} & \mathrm{Fr}\mathbf{X} & \xleftarrow{\mu'} & \mathrm{Fr}^2\mathbf{X} \\ \| & & \downarrow t' & & \downarrow \mathrm{Fr}(t') \\ \mathbf{X} & = & \mathbf{X} & \xleftarrow{t'} & \mathrm{Fr}\mathbf{X} \end{array}$$

Moreover, any $\mathrm{F}r$-algebra $t'\colon \mathrm{F}r\mathbf{X} \to \mathbf{X}$ determines a $\mathrm{P}$-algebra $t = t'p\colon \mathrm{P}\mathbf{X} \to \mathbf{X}$, while a $\mathrm{P}$-algebra $t\colon \mathrm{P}\mathbf{X} \to \mathbf{X}$ induces a $\mathrm{F}r$-algebra $t'\colon \mathrm{F}r\mathbf{X} \to \mathbf{X}$ (with $t = t'p$) iff $t$ is compatible with R.

**3. Pseudo algebras.** Actually, we want to compare the 2-category Fs**Cat** (of categories with fs, functors which preserve them, and natural transformations of such functors) with the 2-category Ps-P-**Alg** of pseudo P-algebras, always understood to be *unitary* (or *normalised*).

According to a general definition (cf. [St], §2), a (unitary) *pseudo* P*-algebra* $(\mathbf{X}, t, \vartheta)$, or *factorisation pseudo algebra*, consists of a category $\mathbf{X}$, a functor $t$ (the *structure*) and a functorial isomorphism $\vartheta$ (*pseudo associativity*), so that

(1) $t\colon \mathbf{X}^\mathbf{2} \to \mathbf{X}, \qquad\qquad t.\eta\mathbf{X} = 1$,

(2) $\vartheta\colon t.\mathrm{P}t \cong t.\mu\mathbf{X}\colon \mathrm{P}^2\mathbf{X} \to \mathbf{X}$,

(3) $\vartheta(\mathrm{P}\eta\mathbf{X}) = 1_t = \vartheta(\eta\mathrm{P}\mathbf{X})\colon t \to t\colon \mathrm{P}\mathbf{X} \to \mathbf{X}$,

(4) $\vartheta(\mathrm{P}\mu\mathbf{X}).t(\mathrm{P}\vartheta) = \vartheta(\mu\mathrm{P}\,\mathbf{X}).\vartheta(\mathrm{P}^2t)\colon t.\mathrm{P}t.\mathrm{P}^2t \to t.\mu\mathbf{X}.\mu\mathrm{P}\mathbf{X}\colon \mathrm{P}^3\mathbf{X} \to \mathbf{X}$,

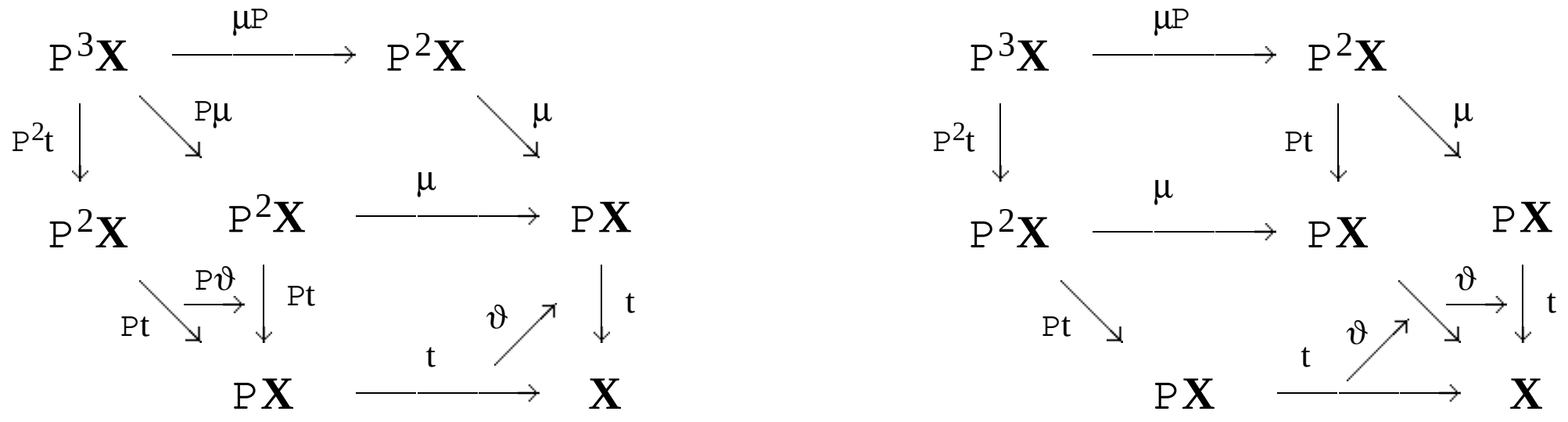

but here (i.e., for P) the conditions (3), (4) follow from the rest (as proved below, 4 (A), (B)).

A *morphism* $(F, \varphi)\colon (\mathbf{X}, t, \vartheta) \to (\mathbf{Y}, t', \vartheta')$ of pseudo P-algebras is a functor $F\colon \mathbf{X} \to \mathbf{Y}$ with a functorial isomorphism $\varphi\colon F.t \to t'.\mathrm{P}F\colon \mathrm{P}\mathbf{X} \to \mathbf{Y}$ satisfying the following coherence conditions (again, the second is redundant for P, cf. 4 (A), (B))

(5) $\varphi.\eta\mathbf{X} = 1_F\colon \mathbf{X} \to \mathbf{Y}$,

(6) $\varphi\mu\mathbf{X}.F\vartheta = \vartheta'\mathrm{P}^2F.t'\mathrm{P}\varphi.\varphi\mathrm{P}t\colon F.t.\mathrm{P}\,t \to t'.\mu\mathbf{Y}.\mathrm{P}^2F\colon \mathrm{P}^2\mathbf{X} \to \mathbf{Y}$.

Finally, a 2-cell $\alpha: (F, \varphi) \to (G, \psi)$ is just a natural transformation $\alpha: F \to G$; it is automatically coherent (cf. 4 (B))

(7) $\psi.\alpha t \; = \; t'P\alpha.\varphi: \; F.t \to t'.PG: \; P\,\mathbf{X} \to \mathbf{Y}.$

Similarly, we have the 2-category Ps-$\mathrm{Fr}$-**Alg** of *pseudo* $\mathrm{Fr}$*-algebras,* or *proper-factorisation pseudo algebras*; these amount to pseudo $P$-algebras $(\mathbf{X}, t, \vartheta)$ where both $t$ and $\vartheta$ are consistent with $R$ (the consistency of $\vartheta$ being redundant, cf. 4 (D).). Again, (3), (4), (6), (7) are redundant.

**4. Theorem (The comparison of factorisation algebras).** (i) (Coppey-Korostenski-Tholen) With respect to the diagonal 2-monad for the endofunctor $P = (-)^2$ of **Cat**, there is a canonical equivalence of categories – described below – between Fs**Cat** and Ps-$P$-**Alg**, which induces a bijection between fs on a category $\mathbf{X}$ and pseudo isomorphism classes of pseudo $P$-algebras on $\mathbf{X}$. In the strict situation, the canonical comparison functor $K_0$: $\mathrm{Fs}_0$**Cat** $\to$ $P$-**Alg**, between strict fs and $P$-algebras, is an isomorphism.

(ii) With respect to the 2-monad of the endofunctor $\mathrm{Fr}$, the previous equivalence induces an equivalence between categories with proper factorisation systems and pseudo $\mathrm{Fr}$-algebras, as well as a bijection between *proper* fs on a category $\mathbf{X}$ and pseudo isomorphism classes of pseudo $\mathrm{Fr}$-algebras on $\mathbf{X}$. The comparison functor $K'_0$: $\mathrm{PFs}_0$**Cat** $\to$ $\mathrm{Fr}$-**Alg**, of proper strict fs, is an isomorphism.

**Proof.** Part (i) is mostly proved in [KT], and we only need to complete a few points.

(A) First, there is a canonical 2-functor $L$: Ps-$P$-**Alg** $\to$ Fs**Cat**. Given a (unitary) pseudo $P$ -algebra $(\mathbf{X}, t, \vartheta)$, every map $x: X' \to X''$ in $\mathbf{X}$ inherits a precise *t-factorisation* through the object $t(\hat{x})$, by letting the functor $t$ act on the canonical factorisation of $\eta\mathbf{X}(x) = (x, 1).(1, x)$ in $\mathbf{X}^2$ (1.3)

(1) $$\begin{array}{ccccc} X' & = & X' & \xrightarrow{x} & X'' \\ {\scriptstyle 1}\downarrow & & \downarrow{\scriptstyle x} & & \downarrow{\scriptstyle 1} \\ X' & \xrightarrow{x} & X'' & = & X'' \end{array} \qquad\qquad X' \xrightarrow{\tau^-(x)} t(\hat{x}) \xrightarrow{\tau^+(x)} X''$$

(2) $\tau^-(\hat{x}) \; = \; t(1, x): \; X' \to t(\hat{x}), \qquad\qquad \tau^+(\hat{x}) \; = \; t(x, 1): \; t(\hat{x}) \to X'',$

$\tau^+(\hat{x}).\tau^-(\hat{x}) \; = \; t((x, x): X' \to X'') \; = \; t.\eta(x) \; = \; x.$

$E$ is defined as the class of $\mathbf{X}$-maps $x$ such that $\tau^+(\hat{x})$ is iso; dually for $M$. This is indeed a fs, as proved in [KT], thm. 4.4, *without* assuming the coherence condition 3.3 (cf. the Note at the end of the paper) *nor* 3.4; the fact that these properties will be obtained in (B), from the backward procedure, shows that they are redundant. (In the strict case, a strict $P$-algebra $t$ gives a strict fs, where $E_0$ contains the maps $x$ such that $\tau^+(\hat{x})$ is an identity, and dually for $M_0$.)

Given a morphism $(F, \varphi): (\mathbf{X}, t, \vartheta) \to (\mathbf{Y}, t', \vartheta')$ of pseudo $P$-algebras, the fact that the functor $F: \mathbf{X} \to \mathbf{Y}$ preserves the associated fs follows from the following diagram, commutative by the naturality of $\varphi: F.t \to t'.PF$ on $(1, x): X' \to \hat{x}$, $(x, 1): \hat{x} \to X''$, $(1, y)$ and $(y, 1)$

(3) 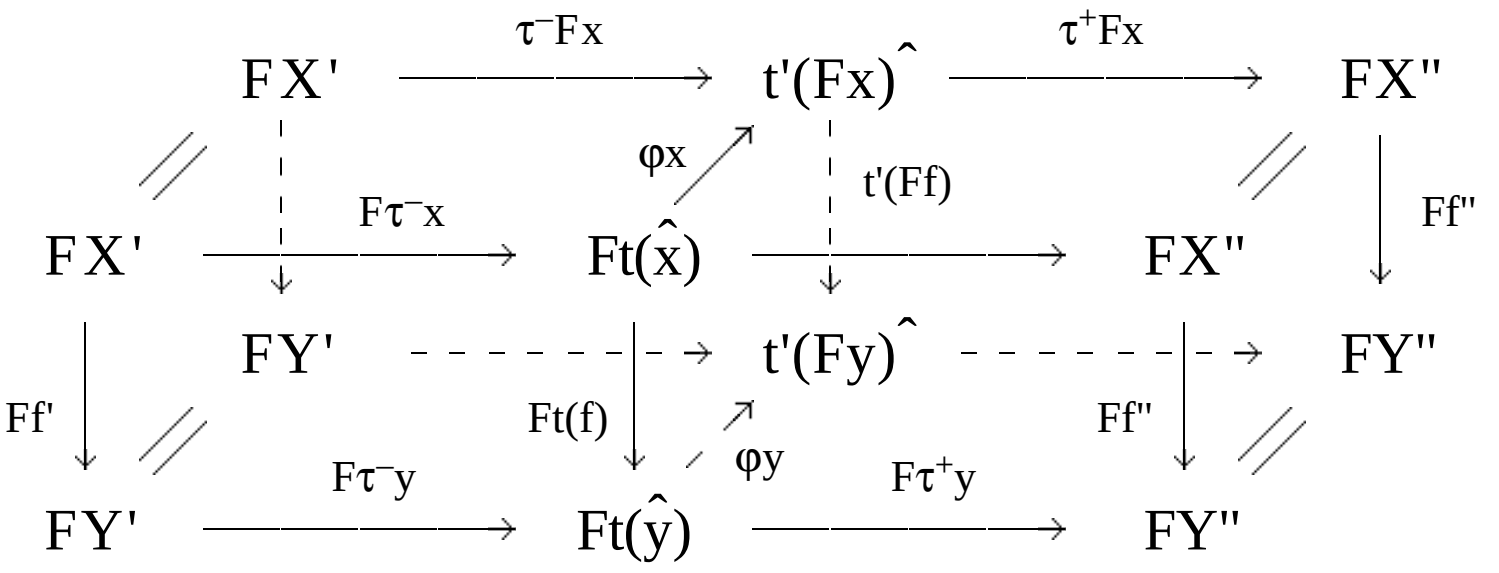

Again, we do *not* need the condition 3.6: any natural iso $\varphi$ such that $\varphi.\eta\mathbf{X} = 1_F$ has this effect.

(B) Conversely, one can construct a 2-functor $K : Fs\mathbf{Cat} \to$ Ps-$P$-**Alg** *depending on choice*. Let $(\mathbf{X}, E, M)$ be a category with fs; for every map $x: X' \to X''$, let us *choose* one structural factorisation $x = \tau^+(x).\tau^-(x): X' \to t(x) \to X''$, respecting all identities: $1 = 1.1$ (We are *not* saying that this choice comes from a strict fs contained in $(E, M)$). By orthogonality, this choice determines *one* functor $t: \mathbf{X}^2 \to \mathbf{X}$ with this action on the objects and such that $\tau^-: \partial^- \to t$, $\tau^+: t \to \partial^+$ are natural transformations ($\partial^-, \partial^+: \mathbf{X}^2 \to \mathbf{X}$ being the domain and codomain functors)

$$(4)\qquad \begin{array}{ccccc} X' & \xrightarrow{\tau^- x} & t(x) & \xrightarrow{\tau^+ x} & X'' \\ {\scriptstyle f'}\downarrow & & \downarrow{\scriptstyle t(f)} & & \downarrow{\scriptstyle f''} \\ Y' & \xrightarrow[\tau^- y]{} & t(y) & \xrightarrow[\tau^+ y]{} & Y'' \end{array} \qquad\qquad f = (f', f''): x \to y.$$

Now $t.\eta(X) = t(1_X) = X$. Moreover, let $t.Pt$ and $t.\mu\mathbf{X}: P^2\mathbf{X} \to \mathbf{X}$ operate on the object $(f', f''): x \to y$ of $P^2\mathbf{X}$, producing $t.Pt(f', f'') = Z'$ and $t.\mu\mathbf{X}(f', f'') = t(\overline{f}) = Z''$

$$(5)\qquad \begin{array}{ccccc} X' & \xrightarrow{x'} & t(x) & \xrightarrow{x''} & X'' \\ & & \downarrow{\scriptstyle z'} & & \\ {\scriptstyle f'}\downarrow & & Z' & & \downarrow{\scriptstyle f''} \\ & & \downarrow{\scriptstyle z''} & & \\ Y' & \xrightarrow[y']{} & t(y) & \xrightarrow[y'']{} & Y'' \end{array} \qquad\qquad \begin{array}{ccccc} X' & & \xrightarrow{x} & & X'' \\ & \searrow{\scriptstyle d'} & & & \\ {\scriptstyle f'}\downarrow & & Z'' & & \downarrow{\scriptstyle f''} \\ & & & \searrow{\scriptstyle d''} & \\ Y' & & \xrightarrow[y]{} & & Y'' \end{array}$$

so that there is precisely one isomorphism $\vartheta(f): Z' \to Z''$ linking the two EM-factorisations we have obtained for the diagonal, $\overline{f} = (y''z'').(z'x') = d''.d'$ (a strict fs would give an identity, for $\vartheta(f)$)

(6) $\vartheta(f): t.Pt(f) \to t.\mu\mathbf{X}(f),\qquad\qquad y''z'' = d''.\vartheta(f),\qquad \vartheta(f).(z'x') = d'.$

The coherence relations for $\vartheta$ do hold: the first (3.3) is obvious; the second (3.4) is concerned with two natural transformations, $\vartheta(P\mu\mathbf{X}).t(P\vartheta)$ and $\vartheta(\mu P\,\mathbf{X}).\vartheta(P^2t)$, that take a commutative cube $\Xi \in \mathrm{Ob}(P^3\mathbf{X})$ to the unique isomorphism linking two precise EM-factorisations of the diagonal arrow of $\Xi$, through $t.Pt.P^2t(\Xi)$ and $t.\mu\mathbf{X}.\mu P\mathbf{X}(\Xi)$, respectively.

By similar arguments, a functor $F: (\mathbf{X}, E, M) \to (\mathbf{Y}, E', M')$ that preserves fs is easily seen to produce a morphism $(F, \varphi): (\mathbf{X}, t, \vartheta) \to (\mathbf{Y}, t', \vartheta')$ of the associated pseudo $P$-algebras. Note that $\varphi: F.t \to t'.PF: P\mathbf{X} \to \mathbf{Y}$ is determined by the choices which give $t$ and $t'$, and does satisfy the coherence condition 3.6, $\varphi\mu\mathbf{X}.F\vartheta = \vartheta'P^2F.t'P\varphi.\varphi Pt$; these two natural transformations take a commutative square $\xi \in \mathrm{Ob}(P^2\mathbf{X})$ to the unique isomorphism linking two precise EM-factorisations of the diagonal arrow of the square, through $F.t.Pt(\xi)$ and $t'.\mu\mathbf{Y}.P^2F(\xi)$. Similarly, a natural transformation $\alpha: F \to G$ satisfies automatically the condition 3.7.

(C) The composite Fs**Cat** $\to$ Ps-P-**Alg** $\to$ Fs**Cat** is the identity. Let $(E, M)$ be a fs on a category **X**, $(t, \vartheta)$ the associated pseudo P-algebra and $(E', M')$ the fs corresponding to the latter. Then $E' = \{x \mid \tau^+(x) \text{ is iso}\}$ plainly coincides with $E$, and $M' = M$.

The other composite, Ps-P-**Alg** $\to$ Fs**Cat** $\to$ Ps-P-**Alg**, is just isomorphic to the identity. It is now sufficient to consider two pseudo P-algebras $(t, \vartheta)$, $(t', \vartheta')$ on **X**, giving the same factorisation system $(E, M)$, and prove that they are pseudo isomorphic, in a unique coherent way. Actually, for each $x: X' \to X''$ in **X** there is one iso $\varphi(x)$ linking the t- and t'-factorisation (both in $(E, M)$)

$$(7) \qquad \begin{array}{ccccc} X' & \xrightarrow{\tau^-x} & t(\hat{x}) & \xrightarrow{\tau^+x} & X'' \\ \| & & \downarrow{\scriptstyle \varphi x} & & \| \\ X' & \xrightarrow{\tau'^-x} & t'(\hat{x}) & \xrightarrow{\tau'^+x} & X'' \end{array}$$

this gives a functorial isomorphism $\varphi: t \to t': \mathrm{P}\mathbf{X} \to \mathbf{X}$ such that $(1_{\mathbf{X}}, \varphi): (\mathbf{X}, t, \vartheta) \to (\mathbf{X}, t', \vartheta')$ is a pseudo isomorphism of algebras.

(D) For Part (ii), we only need now to prove that, in the previous transformations, pseudo Fr-algebras (i.e., pseudo P-algebras consistent with the Freyd congruence R) correspond to proper fs.

First, the consistency of $t: \mathbf{X}^2 \to \mathbf{X}$ with R is sufficient to give an epi-mono factorisation system. Take, for instance, $m \in M$ (so that $u = \tau^-(m)$ is iso) and $mf_1 = h = mf_2$ in the left-hand diagram below; then, the naturality of the transformation $\tau^-: \partial^- \to t$ on the R-equivalent maps $(f_i, h): X' \to m$ of $\mathbf{X}^2$ gives $uf_1 = t(f_1, h) = t(f_2, h) = uf_2$ and $f_1 = f_2$

$$(8) \qquad \begin{array}{ccc} X' & = & X' \\ {\scriptstyle f_i}\downarrow & & \downarrow{\scriptstyle h} \\ X & \xrightarrow{m} & Y \end{array} \qquad\qquad \begin{array}{ccc} X' & = & X' \\ {\scriptstyle f_i}\downarrow & & \downarrow{\scriptstyle t(f_i,\, h)} \\ X & \xrightarrow{u} & t(m) \end{array}$$

Finally, if $(E, M)$ is epi-mono, then $t(f)$ in (3) only depends on the diagonal $\overline{f}$ of $f = (f', f'')$: $x \to y$ in $\mathbf{X}^2$, and similarly for $\vartheta(f)$ in (5). Therefore they induce a functor $t'$: Fr**X** $\to$ **X** and a functorial iso $\vartheta'$: t'.Fr(t') $\to$ t'.μ'**X**, which form a pseudo Fr-algebra.

**5. Remarks.** A crucial tool for the proof of point (A), above, is the structure of $\mathrm{P}\mathbf{X} = \mathbf{X}^2$ as a "path functor" (representing natural transformations): it forms a *cubical comonad* [G1, G2], well linked to the previous monad structure. This interplay already arises in the exponent category **2** – a *comonoid* and a *lattice* (more precisely, a *cubical monoid* [G1]) – and was exploited in this form in [KT], Section 1.

The cubical comonad structure, relevant for formal homotopy theory [G2], has one *degeneracy* $\eta$: $1 \to \mathrm{P}$ (the previous unit), two *faces* or *co-units* $\partial^{\pm}: \mathrm{P} \to 1$ (domain and codomain) and two *connections* or *co-operations* $g^{\pm}: \mathrm{P} \to \mathrm{P}^2$

$$(1) \qquad \begin{array}{ccc} X' & \xrightarrow{x} & X'' \\ {\scriptstyle x}\downarrow & g^-(\hat{x}) & \| \\ X'' & = & X'' \end{array} \qquad\qquad \begin{array}{ccc} X' & = & X' \\ \| & g^+(\hat{x}) & \downarrow{\scriptstyle x} \\ X' & \xrightarrow[x]{} & X'' \end{array}$$

(The connections have appeared above in the canonical factorisation $\eta\mathbf{X}(x) = g^{-}\mathbf{X}(\hat{x}).g^{+}\mathbf{X}(\hat{x})$; the natural transformations $\tau^{-}$, $\tau^{+}$ can thus be obtained as $\tau^{-} = \mathrm{P}t.g^{+}\mathbf{X}$, $\tau^{+} = \mathrm{P}t.g^{-}\mathbf{X}$.)

A cubical comonad satisfies axioms [G1, G2] essentially saying that $\partial^{\varepsilon}$ $(\varepsilon = \pm)$ is a co-unit for the corresponding connection $g^{\varepsilon}$ and co-absorbant for the other, while $\eta$ makes everything degenerate; moreover, the connections are co-associative. Here the two structures, monad and cubical comonad, are linked by some equations (after the coincidence of the monad-unit with the degeneracy; the last formula is actually a consequence of the co-associativity of connections):

(2) $\partial^{\varepsilon}\mu = \partial^{\varepsilon}.\mathrm{P}\partial^{\varepsilon} = \partial^{\varepsilon}.\partial^{\varepsilon}\mathrm{P}$, $\quad\mu g^{\varepsilon} = 1_{\mathrm{PX}}$,

$\mathrm{P}\mu.g^{\varepsilon}\mathrm{P}.g^{\varepsilon'} = \eta\mathrm{P}$, $\quad\mathrm{P}\mu.g^{\varepsilon}\mathrm{P}.g^{\varepsilon} = \mathrm{P}\mu.\mathrm{P}g^{\varepsilon}.g^{\varepsilon} = g^{\varepsilon}$ $\quad(\varepsilon \neq \varepsilon')$.

A natural question arises – if the previous arguments have a non-trivial rebound in the usual range of homotopy, the category **Top** of topological spaces. Replace the categorical interval **2** with the topological one, $I = [0, 1]$, which is, again, a diagonal comonoid and a lattice (and an exponentiable object); thus, the path functor $PX = X^{I}$ is a monad and a cubical comonad, consistently as above. But here, the interest of (pseudo?) P-algebras is not clear (once we have excluded the trivial, "universal" ones: for a fixed $a \in I$, every space X has an obvious strict structure, $ev_{a}: PX \to X$; in the same way as each category $\mathbf{X}$ has two trivial P-algebras, $\partial^{\pm}: \mathrm{P}\mathbf{X} \to \mathbf{X}$, and two trivial fs). On the other hand, one can readily note that the Kleisli category of P has for morphisms the homotopies, with "diagonal" horizontal composition: $(\beta\circ\alpha)(x; t) = \beta(\alpha(x; t); t)$, for $t \in I$.